\documentclass[11pt]{article}

\usepackage[T1]{fontenc}
\usepackage{lmodern}
\usepackage[a4paper,margin=1.15in]{geometry}
\usepackage{amsmath,amssymb,amsthm,mathtools}
\usepackage{enumitem}
\usepackage[expansion=false]{microtype}
\usepackage{cite}
\usepackage[hidelinks]{hyperref}
\hypersetup{
  pdftitle={Smith's slope problem and finite variations of analytic capacity},
  pdfauthor={Deguang Zhong}
}

\newtheorem{theorem}{Theorem}[section]
\newtheorem{lemma}[theorem]{Lemma}
\newtheorem{corollary}[theorem]{Corollary}
\newtheorem{proposition}[theorem]{Proposition}
\newtheorem{question}[theorem]{Question}
\theoremstyle{remark}
\newtheorem{remark}[theorem]{Remark}

\theoremstyle{plain}
\newtheorem{maintheorem}{Theorem}

\newcommand{\C}{\mathbb C}
\newcommand{\D}{\mathbb D}
\newcommand{\R}{\mathbb R}
\newcommand{\Chat}{\widehat{\mathbb C}}
\newcommand{\Hm}{\mathcal H}
\newcommand{\dist}{\operatorname{dist}}
\newcommand{\diam}{\operatorname{diam}}

\newcommand{\dd}{\,d}
\newcommand{\dbar}{\overline\partial}
\newcommand{\Fill}{\operatorname{fill}}
\newcommand{\caplog}{\operatorname{cap}_{\log}}
\newcommand{\Fcal}{\mathcal C}

\title{\bfseries Smith's slope problem and finite variations of analytic capacity}

\author{%
Deguang Zhong\\[0.7em]
\small Institute of Applied Mathematics, Shenzhen Polytechnic University, Shenzhen 518055, China\\ \texttt{huachengzhon@163.com}
}

\date{}

\begin{document}

\maketitle

\begin{abstract}
Smith's slope measures the first-order increase in analytic capacity when a small disk is adjoined to a compact plane set. In this paper, we first prove that
this slope is at most one at every finite exterior point and characterize equality for finitely connected exteriors. The proof represents the square of the slope for exteriors with analytic boundary as the Fredholm
determinant of a positive contraction. Secondly, we show that the slope bound is equivalent to a comparison of the capacity increments of two nested
continua in the presence of an arbitrary compact set. This gives exact subadditivity when one of the two sets is connected and strong subadditivity when one is contained in a line. For a compact set of $\sigma$-finite length whose rectifiable part has finite length, we bound the capacity lost on passing to a compact subset by one quarter of the rectifiable length removed. Finally, we bound the slope in
terms of the weighted area of the filled compact set and show that the first-order increment for a compact perturbation of arbitrary shape equals its analytic capacity times the slope at the point of insertion.
\end{abstract}

\medskip
\noindent\textbf{Keywords.} Analytic capacity,  Ahlfors function, Garabedian function, subadditivity, radial slits,  boundary variation, Cauchy capacity, rectifiability.

\medskip
\noindent\textbf{2020 Mathematics Subject Classification.}
Primary 30C85, 30E20. Secondary 30C35, 30C62, 28A75.

\section{Introduction}\label{sec1}

\subsection{From subadditivity to Smith's slope}
Suppose that $E\subseteq\mathbb{C}$ is a compact set and $\Omega_E$ is the component of $\Chat\setminus E$ containing infinity, then the definition of analytic capacity $\gamma,$ which introduced by Ahlfors \cite{A1947}
to characterize the removable sets for bounded holomorphic functions in Painlev\'e's problem \cite{Pa1888}, is defined by 
\begin{equation}\label{eq1}
 \gamma(E)=\sup\{|f'(\infty)|: f\in H^\infty(\Omega_E),\
 \|f\|_\infty\le1,\ f(\infty)=0\}, f'(\infty)=\lim_{z\to\infty}zf(z).
\end{equation}
Here $\Chat$ denotes the Riemann sphere and $H^\infty(\Omega_E)$ is the space of bounded holomorphic functions on $\Omega_E$. The normalized extremizer $f_E$, with $f_E'(\infty)=\gamma(E)$, is 
called the Ahlfors
function. 

For analytic capacity, there is a \textbf{\textit{semiadditivity problem}} arose  in Vitushkin's work on rational approximation \cite[Chapter~I, Section~1]{V1967}; see also Davie \cite[Section~5]{Da1972}.
It asks whether there is a universal constant $M<\infty$ such that
\begin{equation}\label{a}
 \gamma(E\cup F)\le M\bigl(\gamma(E)+\gamma(F)\bigr)
\end{equation}
 holds for all compact $E,F\subset\C$. Tolsa \cite[Theorem~1.4]{T2003} proved the stronger, countable form of semiadditivity.  The exact \textbf{\textit{subadditivity problem}}  for analytic capacity, explicitly posed in Smith's 1974 thesis \cite[Section~19, question~(2), p.~116]{ST1974} (see also Tolsa \cite[Section~5, Problem~1, p.~1521]{T2006}), asks whether one can take $M=1$ in \eqref{a}. Namely, whether the following inequality  
\begin{equation}\tag{SA}\label{eq:SA}
 \gamma(E\cup F)\le\gamma(E)+\gamma(F)
\end{equation}
holds for all compact $E,F\subset\C$. For partial results, the reduction to equal radius disks, and numerical evidence, see \cite{Su1984,YR2013,NGV2025} and the survey \cite[Section~4]{YS2018}.
Smith studied this question by adjoining a single small disk. For a finite point $z\in\Omega_E$, let $\overline D(z,s)$ denote the closed disk of radius $s$ centered at $z$. Smith \cite[Theorems~2.2 and~3.1]{S1974} proved the existence of the \textbf{\textit{slope}}  
\begin{equation}\label{eq2}
 a_E(z)=\lim_{s\downarrow0}
 \frac{\gamma(E\cup\overline D(z,s))-\gamma(E)}s .
\end{equation}
Since $\gamma(\overline D(z,s))=s$, applying \eqref{eq:SA} with $F=\overline D(z,s)$ would bound the difference quotient in
\eqref{eq2} by one for every $s>0$. Smith made this implication explicit in \cite[p.~300]{S1974} and listed the infinitesimal bound as question~(6) in his dissertation \cite[Section~19, p.~116]{ST1974}.

\begin{question}[Smith]\label{que1}
Does the bound
\begin{equation*}\tag{S}\label{eqS}
   a_E(z)\le 1
\end{equation*}
hold for every compact $E\subset\C$ and every finite $z\in\Omega_E$?
\end{question}

For connected $E$, Smith derived \hyperref[eqS]{(S)} from the classical distortion theorem for univalent functions \cite[p.~300]{S1974}.
He also obtained the bound for compact subsets of a line and for sets with at most two components \cite[Corollaries~16.3 and~18.3]{ST1974}. Our first result answers Question~\ref{que1} for arbitrary compact sets.
\begin{maintheorem}\label{thm1}
For every compact $E\subset\C$ and every finite exterior point
$z\in\Omega_E$, one has $0<a_E(z)\le1$. The universal
constant $1$ is attained.
\end{maintheorem}
Theorem \ref{thm1} also answers Smith's question about a universal upper bound for the slope \cite[Section~19, question~(7), p.~116]{ST1974}.
Actually, the sharp constant in Theorem \ref{thm1} is attained, for example, when $E$ is a segment and $z$ lies on its supporting line. To prove the upper bound, we use Smith's expression for the slope in terms of the Ahlfors and Garabedian functions. On an exterior bounded by analytic curves,
Garabedian's dual extremal problem \cite{G1949} gives a boundary identity for the extremizers. We use it to express $a_E(z)^2$ at an exterior point as a Fredholm determinant. The operator in that formula is the
compression of an orthogonal projection to the graph of a bounded operator, so its determinant is at most one. Decreasing approximation then gives the bound for arbitrary compact sets.

Murai  \cite[Theorem~1 and p.~199]{M1992} related the slope to the first variation of capacity along an analytic slit and used Rengel's derivative inequality to reduce
the upper bound to radial slits; see also \cite{R1933}. In Section~\ref{sec3} we express the change in slope under radial-slit
uniformization in terms of a Dirichlet integral and a weighted area integral. This identity and the projection formula characterize
equality in finite connectivity. Equality holds precisely when the non-singleton components of the complement of $\Omega_E$ are straight segments on a single line through $z$; see Corollary~\ref{cor1}.

\subsection{From the slope bound to finite increments}
How much of the finite inequality \eqref{eq:SA} follows from \hyperref[eqS]{(S)}? For two continua, exact subadditivity was established by Smith \cite[Theorem~18.4]{ST1974}. For intersecting continua it also follows from Pommerenke's inequality  \cite[Theorem~2]{P1968} for logarithmic capacity, since the union is connected and the two capacities agree on continua. In 1984, Suita \cite[Theorem~2]{Su1984} gave a conformal mapping proof for two disjoint continua using a collinear slit extremal domain. Here a \emph{continuum} is a nonempty compact connected set, possibly a singleton. The following result compares the cost of enlarging a continuum before and after adjoining an arbitrary compact set $E$.
\begin{maintheorem}\label{thm2}
The following assertions are equivalent.
\begin{enumerate}[label=\textup{(\roman*)},leftmargin=*]
\item The slope bound \hyperref[eqS]{(S)} holds.
\item For every compact $E$ and every closed disk $D$,
 $\gamma(E\cup D)\le\gamma(E)+\gamma(D)$.
\item For every compact $E$ and every pair of continua $C_0\subset C_1$,
 \begin{equation}\label{eq3}
   \gamma(E\cup C_1)-\gamma(E\cup C_0)
      \le\gamma(C_1)-\gamma(C_0).
 \end{equation}
\end{enumerate}
\end{maintheorem}
The main step in the proof of Theorem~\ref{thm2} is the implication \textup{(i)}$\Rightarrow$\textup{(iii)}.  To do this, we first establish Proposition~\ref{prop3}, where we use a normalized Joukowski
coordinate to compare $a_B(z)$ with $a_C(z)$ when $C\subset B$ is a continuum. Then, in the separated analytic case, we apply this comparison with $B=E\cup C_t$ and $C=C_t$, where $C_t$ expands from $C_0$ to $C_1$.
The two slopes vanish on the moving boundary, so their pointwise order gives the same order for their exterior normal derivatives. The first variation formula in Lemma~\ref{lem4} thus gives
\[
   \frac{\dd}{\dd t}\gamma(E\cup C_t)
      \le \frac{\dd}{\dd t}\gamma(C_t).
\]
The mean value theorem gives \eqref{eq3} in this case. In Section~\ref{sec4}, we  extend this comparison by approximation and treat intersections separately. Taking $C_0$ to be a point of $C_1$ gives \eqref{eq:SA} whenever one summand is a continuum.

A further consequence concerns a possibly disconnected compact subset $F$ of a line. By Pommerenke's formula \cite{P1960}, $\gamma(F)=\Hm^1(F)/4$, where $\Hm^1$ stands for one-dimensional Hausdorff measure normalized to agree with length on intervals. Applying Theorem~\ref{thm2} to the individual intervals in a cover of $F\setminus E$ gives the strong subadditivity inequality
\[
 \gamma(E\cup F)+\gamma(E\cap F)\le\gamma(E)+\gamma(F)
\]
for every compact $E$; see also in Corollary~\ref{cor3}.

Using Mel'nikov's discretization \cite{Me1995}, Younsi and Ransford \cite[Theorem~7.1]{YR2013} reduced \eqref{eq:SA} to disjoint compact sets $E$ and $F$ that are finite unions of pairwise disjoint closed disks, all of the same radius $r$; see \cite{NGV2025} for numerical examples. If $F$ consists of $k$ such disks, adjoining them one at a time by Theorem~\ref{thm2} gives $\gamma(E\cup F)-\gamma(E)\le kr$. The bound required in \eqref{eq:SA} is $\gamma(F)$.

In the following result, the factor $1/4$  controls how much capacity can be lost when part of a compact set is removed, provided its rectifiable part has finite length. The geometric
condition is the removability of purely unrectifiable compact sets of finite length (David \cite{D1998}) and of $\sigma$-finite length (Tolsa \cite{T2003}). Here, we use the rectifiability conventions of
Section~\ref{sec2}.

\begin{maintheorem}\label{thm3}
Let $K=R\cup P$ be a compact set of $\sigma$-finite $\Hm^1$ measure, where $R,P$ are Borel, $R$ is countably $1$-rectifiable with $\Hm^1(R)<\infty$, and $P$ is purely $1$-unrectifiable. For every compact $K_0\subset K$,
\begin{equation}\label{eq4}
 0\le\gamma(K)-\gamma(K_0)\le\tfrac14\Hm^1(R\setminus K_0).
\end{equation}
The constant is sharp. Consequently, if $K_n$ is an increasing sequence of compact sets with $\bigcup_n K_n=K$, then $\gamma(K_n)\uparrow\gamma(K)$.
Moreover, $\gamma_c(K)=\gamma(K)$.
\end{maintheorem}

The continuity of analytic capacity from above for decreasing compact sets follows from a normal family argument. Younsi \cite[Theorem~3.2]{Y2018} proved convergence of the Ahlfors functions, and hence of their derivatives at infinity, under Carath\'eodory convergence of nondegenerate exteriors of the same finite connectivity. He also asked whether $\gamma(K_n)\uparrow\gamma(K)$ for every increasing sequence of compact sets whose union $K$ is compact \cite[Problem~6.11]{YS2018}; see the concluding remark of \cite{Y2018}. Theorem~\ref{thm3} gives a quantitative answer under its length hypotheses, without any restriction on connectivity. The total length of $K$ may be infinite, and \eqref{eq4} bounds the loss with sharp constant $1/4$ using only the length removed from the rectifiable part.
The same estimate gives $\gamma_c(K)=\gamma(K)$ by approximation with compact subsets of the rectifiable part, where Havinson's theorem for sets of finite Painlev\'e length applies \cite{H1961,Y2018}.

In the rest of this paper, Section~\ref{sec2} discusses the capacity conventions, the dual extremal formulas, and the conformal covariance of the slope. In Section~\ref{sec3}, we prove the slope bound by the graph determinant and use radial-slit coordinates to characterize equality in finite connectivity. Section~\ref{sec4} integrates the pointwise comparison to obtain the finite increment inequality and its subadditivity consequences. In Section~\ref{sec5},  the sharp relative length estimate and the continuity and Cauchy-capacity conclusions are obtained. Section~\ref{sec6} derives area bounds for the slope and the first-order expansion for small compact perturbations.

\section{Capacity, duality, and conformal coordinates}\label{sec2}

\subsection{Definitions and conventions}
We use the conventions of \cite{Ga1972,T2014}. The \emph{filled compact set} is defined by
\[
\Fill(E):=\widehat{\mathbb C}\setminus\Omega_E,
\]
where $\Omega_E$ denotes the component of $\widehat{\mathbb C}\setminus E$ containing $\infty$. Equivalently, $\Fill(E)$ is obtained from $E$ by adjoining all bounded components of its complement.
The inclusions
\[
 E\cup F\subset\Fill(E)\cup\Fill(F)\subset\Fill(E\cup F)
\]
show that filling either summand preserves the exterior of the union. The definition of $\gamma$ gives the same value if bounded holomorphic functions on the full complement are used, since one
can extend an exterior test by zero on the other components.

Monotonicity and continuity under decreasing compact approximation follow from \eqref{eq1} and normal families \cite{A1947,Ga1972,T2014}. To approximate a nonempty filled compact from above, take nested regular neighborhoods bounded by finitely many disjoint analytic Jordan curves and fill them. The intersection is the original compact. Indeed, a point of its
exterior can be joined to infinity by a path in the exterior which avoids all sufficiently late neighborhoods. For a nondegenerate continuum one can use level sets of the exterior Riemann map.

For a nonempty compact set $E$, the logarithmic capacity \cite[Definition~1.8]{Sa2010} is
\[
 \caplog(E)=\exp\!\left(-\inf_{\mu}
     \iint_{E\times E}\log\frac{1}{|u-v|}\,\dd\mu(u)\dd\mu(v)\right),
\]
where $\mu$ ranges over Borel probability measures on $E$ and $\exp(-\infty)=0$. By the Fekete--Szeg\H{o} theorem \cite[Theorem~1.9]{Sa2010},
it equals the transfinite diameter, the limit as $m\to\infty$ of
\[
 \sup_{w_1,\ldots,w_m\in E}
       \left(\prod_{i<j}|w_i-w_j|\right)^{2/(m(m-1))}.
\]
On a nondegenerate continuum, $\gamma=\caplog$ is the coefficient of $1/z$ of the normalized exterior map to $\D$, with value zero at infinity and positive derivative there. Both capacities vanish on a singleton \cite{Ga1972,T2014}.

Following \cite{T2003,T2014,Y2018,YS2018}, for a finite complex Borel measure $\mu$ with compact support we write $\Fcal\mu(z)=\int (z-w)^{-1}\dd\mu(w)$ for its \emph{Cauchy transform}, defined almost everywhere and holomorphic off the support. The \emph{Cauchy capacity} $\gamma_c(E)$ is the supremum of $|\mu(E)|$ over such measures supported on $E$ for which $|\Fcal\mu|\le1$ on $\C\setminus E$. The \emph{positive Cauchy capacity} $\gamma_+(E)$ takes the supremum of $\mu(E)$ over positive measures supported on $E$ with $\|\Fcal\mu\|_{L^\infty(\C,\dd A)}\le1$, where $\dd A$ denotes
planar area measure. Thus $\gamma_+\le\gamma_c\le\gamma$.

In Havinson's criterion \cite{H1961,Y2018}, a compact set $E$ has \emph{finite Painlev\'e length} if, in any neighborhood of $E$, one can enclose $E$ by a finite cycle of analytic Jordan curves
disjoint from $E$ with uniformly bounded total length, having winding number $1$ about each point of $E$. Finite $\Hm^1(E)$ implies this property. Havinson proved $\gamma_c(E)=\gamma(E)$
for sets of finite Painlev\'e length. In Smith's terminology \cite[Section~19]{ST1974},  \textbf{\textit{strong subadditivity}}  is
$$\gamma(E\cup F)+\gamma(E\cap F)\le\gamma(E)+\gamma(F).$$

We follow Federer \cite[Section~3.2.14]{Fe1969} and Mattila \cite[Definition~15.3 and Lemma~15.4]{Ma1995} for the following rectifiability conventions. A Borel set $R\subset\C$ is
\emph{countably $1$-rectifiable} if there are sets $A_j\subset\R$ and Lipschitz maps $g_j:A_j\to\C$ such that
\[
   \Hm^1\!\left(R\setminus\bigcup_{j=1}^{\infty}g_j(A_j)\right)=0 .
\]
A Borel set $P\subset\C$ is \emph{purely $1$-unrectifiable} if $\Hm^1(P\cap\Gamma)=0$ for every rectifiable curve $\Gamma$. For the finite-length rectifiable part in Theorem~\ref{thm3}, the
covering sets may be taken to be embedded $C^1$ subarcs up to an $\Hm^1$-null set \cite[Theorem~15.21]{Ma1995}.

\subsection{The Garabedian extremizers}

Assume temporarily that the exterior of a filled compact set has finitely many analytic Jordan boundary curves. We call such
an exterior analytic. Following \cite{G1949,H1961,YR2013}, normalize the Garabedian function by $\psi_E(\infty)=1/(2\pi i)$ and its square
root by
\[
  q_E^2=2\pi i\psi_E,\qquad q_E(\infty)=1,\qquad h_E=f_Eq_E .
\]
Write $I_E(g)=(2\pi)^{-1}\int_{\partial\Omega_E}|g|^2\dd s$. When the boundary is oriented counterclockwise about the compact
components, write $T$ for its unit tangent. The Garabedian boundary identity is $h_ET=i\overline{q_E}$.

\begin{lemma}\label{lem1}
For an analytic exterior, with holomorphic test functions regular up to the boundary,
\begin{equation}\label{eq5}
 \gamma(E)=\min_{g(\infty)=1} I_E(g)
 =\max_{h(\infty)=0}\bigl(2\Re h'(\infty)-I_E(h)\bigr).
\end{equation}
The extremizers are $q_E,h_E$, respectively. In particular $I_E(q_E)=I_E(h_E)=\gamma(E)$, and
\[
  q_E(z)=1+\frac{b_E}{z^2}+O(z^{-3}),\;\;\;
  h_E(z)=\frac{\gamma(E)}z+O(z^{-2}).
\]
For arbitrary compact $E$, decreasing analytic outer approximation defines $q_E,h_E$ by locally uniform limits, with the same normalizations at infinity. The functions $q_E$ are zero-free.
One has $a_E>0$ and $a_E=|q_E|^2-|h_E|^2$ and, for small $s>0$,
\begin{equation}\label{eq:uniform-expansion}
 \gamma(E\cup\overline D(z,s))=\gamma(E)+s \cdot a_E(z)+O(s^2).
\end{equation}
The error is uniform for a fixed positive lower bound on $\dist(z,E)$ and a fixed upper bound on $\gamma(E)$.
\end{lemma}
\begin{proof}
For $E=\varnothing$ take $q_E=1$, $h_E=0$ and $a_E=1$. For analytic exteriors, the first extremal identity and the boundary relation follow from Garabedian's duality \cite{G1949}. For admissible trials $G,H$ in the minimum and maximum problems, put $J_TG=i\overline{GT}$ on the boundary. This is an anti-linear isometry, with $J_Tq_E=h_E$.
The residue formula applied to $GH$ gives
\begin{equation}\label{eq:trial-gap}
 I_E(G)-\bigl(2\Re H'(\infty)-I_E(H)\bigr)=I_E(H-J_TG).
\end{equation}
Taking $G=q_E$ proves the maximum formula, with equality at $h_E$, as in \cite[Theorem~3.2]{YR2013}. The orthogonality
$\langle q_E,h_E\rangle_{I_E}=0$ follows by varying the first minimum along $h_E$. Equivalently, by the boundary identity, $\int_{\partial\Omega_E}q_E^2\dd z=0$, which removes the
$z^{-1}$ coefficient of $q_E$. The residue formula for $q_E^2$ and $I_E(q_E)=\gamma(E)$ give
\begin{equation}\label{eq6}
 |q_E(z)|^2
    \le 1+\frac{\gamma(E)}{\dist(z,E)}\;\;\;(z\in\Omega_E\setminus\{\infty\})
\end{equation}
for analytic exteriors. Since $|h_E|\le |q_E|$, Cauchy's estimates then control both functions and their first derivatives near $z$ in terms of the bounds in the statement.

Write $q=q_E$, $h=h_E$, $q_0=q(z)$ and $h_0=h(z)$.
For the exterior obtained by adding $\overline D(z,s)$, choose
\[
 G(w)=q(w)+sA\frac{q(w)}{w-z},\;
 H(w)=h(w)+sB\frac{q(w)}{w-z},\;
 A=\frac{\overline{h_0}}{q_0},\;B=\frac{\overline{q_0}}{q_0}.
\]
If $q_0=0$ then $h_0=0$, and we use $A=B=0$. In either case $|A|,|B|\le1$. On the old boundary $I_E(H-J_TG)=O(s^2)$ because there $|w-z|\ge\dist(z,E)$ and
$I_E(q)=\gamma(E)$. On the new circle $w=z+se^{it}$, $J_TG=e^{-it}\overline G$, and
\[
 G=q_0+\overline{h_0}e^{-it}+O(s),\;
 H=h_0+\overline{q_0}e^{-it}+O(s).
\]
The circle contributes $O(s^3)$ to the gap in \eqref{eq:trial-gap}. To expand the upper trial value, use
\[
 \frac1{2\pi}\int_{\partial\Omega_E}
             \frac{|q(w)|^2}{w-z}\dd s(w)
   =\frac1{2\pi i}\int_{\partial\Omega_E}\frac{q(w)h(w)}{w-z}\dd w
   =-q_0h_0 .
\]
The negative sign comes from the pole at $z$ in the exterior. The old boundary contributes $\gamma(E)-2s|h_0|^2+O(s^2)$, and the new circle contributes $s(|q_0|^2+|h_0|^2)+O(s^2)$.
This proves the uniform expansion for analytic exteriors. The constants in the old-boundary estimates are controlled by $\dist(z,E)^{-1}$ and $I_E(q)=I_E(h)=\gamma(E)$.
On the new circle the local Cauchy estimates from \eqref{eq6} give the same dependence. Thus the error bound and the allowed range of $s$ are uniform when the distance has a fixed positive lower bound and the capacity has a fixed upper bound.

Now take filled analytic outer approximants $E_n\downarrow\Fill(E)$. Filling $E$ preserves its exterior and the capacities of both $E$ and $E\cup\overline D(z,s)$.
For fixed small $s>0$, both $\gamma(E_n)$ and $\gamma(E_n\cup\overline D(z,s))$ converge by continuity from above. Their difference quotient differs from
$a_n(z)=|q_{E_n}(z)|^2-|h_{E_n}(z)|^2$ by at most $Cs$, uniformly for all sufficiently large $n$. The limsup and liminf of $a_n(z)$ therefore differ by at most $2Cs$.
Letting $s\downarrow0$ shows that $a_n(z)$ has a limit $a(z)$. For each fixed sufficiently small $s>0$, we may now let $n\to\infty$ in the same estimate to obtain
\[
 \bigl|\gamma(E\cup\overline D(z,s))-\gamma(E)-s a(z)\bigr|\le Cs^2 .
\]
Dividing by $s$ and letting $s\downarrow0$ gives $a(z)=a_E(z)$ and proves \eqref{eq:uniform-expansion}.

By \eqref{eq6}, the $q_{E_n}$ form a normal family on compact subsets of $\Omega_E$. The normalized Ahlfors functions converge to $f_E$: a normal limit is admissible, has derivative $\gamma(E)$ at infinity, and is unique. Since $|f_E|<1$ at each point of $\Omega_E$, any limit $q$ satisfies $|q|^2=a_E/(1-|f_E|^2)$ and $q(\infty)=1$. A holomorphic function with this prescribed modulus and normalization is unique. This proves the stated convergence of $q_{E_n}$, $h_{E_n}$ and the slopes, and identifies $\psi_E=q_E^2/(2\pi i)$ with Smith's normalization. Taking limits in \eqref{eq6} preserves that estimate. The $q_{E_n}$ are zero-free by the classical analytic logarithm for the Garabedian function \cite[Theorem~2.2]{YR2013}. Hurwitz's theorem and $q_E(\infty)=1$ imply that $q_E$ is zero-free as well. Consequently $a_E=|q_E|^2(1-|f_E|^2)>0$ on $\Omega_E$. If $\gamma(E)=0$, the estimate and the normalization give $q_E=1$ and $h_E=0$.
\end{proof}

The same uniform expansion also gives convergence for any decreasing sequence of compacts $E_n\downarrow E$, as in \cite[Theorem~3.1]{S1974}. Indeed, for each $z\in\Omega_E$ and fixed small $s$, the difference
quotients converge by continuity from above. The error after replacing the quotients by the slopes is bounded by $Cs$ uniformly in $n$. Let $s\downarrow0$ and use the preceding normal-family argument to
obtain locally uniform convergence of $q_{E_n},h_{E_n}$ and $a_{E_n}$.

\begin{lemma}\label{lem2}
Let $\Phi:\Omega_E\to\Omega_{E'}$ be a conformal bijection with $\Phi(z)=z+O(1)$ near infinity. Then $\gamma(E)=\gamma(E')$, and $a_E(z)=|\Phi'(z)|a_{E'}(\Phi(z))$ at finite exterior points.
\end{lemma}
\begin{proof}
Composition preserves the derivative at infinity under this normalization and gives $f_E=f_{E'}\circ\Phi$. For sufficiently small $s>0$, the image of $\overline D(z,s)$ is sandwiched between concentric disks centered
at $\Phi(z)$ of radii $(|\Phi'(z)|+o(1))s$. Removing the disk before or after mapping gives conformally equivalent exteriors. Apply monotonicity and \eqref{eq2}. In particular the slope transformation
law and Lemma~\ref{lem1} give
\begin{equation}\label{eq:q-covariance}
 q_E^2=(q_{E'}\circ\Phi)^2\Phi'.
\end{equation}
The quotient of the two sides has constant modulus $1$, by $a_E=|q_E|^2(1-|f_E|^2)$, and its value at infinity is $1$.
\end{proof}

\section{The slope bound}\label{sec3}
\subsection{Reduction to radial slits}

For a finitely connected exterior, let $S_E$ be the finite set of points that are singleton components of $\Fill(E)$ and put
$$E^\sharp=\Fill(E)\setminus S_E,\;\;\; \Omega_E^\sharp=\Omega_E\cup S_E.$$
The set $E^\sharp$ is compact, as it is the union of the remaining finitely many components.
Filling $E$ preserves capacity, and removing $S_E$ from the filled set preserves it by removability of isolated points. The same holds after adjoining a small disk centered at a finite point $z\in\Omega_E$.
Indeed, each point of $S_E$ still exposed in the exterior of the enlarged compact is an isolated removable singularity. Thus $a_E(z)=a_{E^\sharp}(z)$.
The classical two-point radial-slit uniformization \cite[Chapter~VII, p.~335]{N1952} gives a conformal bijection $\Phi:\Omega_E^\sharp\to\Chat\setminus J$ with
$\Phi(z)=0$ and $\Phi(\zeta)=\zeta+O(1)$ at infinity. Here $J$ is a finite union of disjoint radial segments $[a,b]\tau$, with $0<a<b$ and $|\tau|=1$.
If $E^\sharp=\varnothing$, take $\Phi(\zeta)=\zeta-z$ and $J=\varnothing$. The following identity supplies the derivative estimate in this normalization.

\begin{lemma}\label{lem3}
Let $F=\Phi^{-1}:\Chat\setminus J\to\Omega_E^\sharp$ and write $u(w)=\log|F(w)-z|-\log|w|$, with its continuous value at zero. Then
\begin{equation}\label{eq:energy}
 2\pi\log|F'(0)|
   =\int_{\C\setminus J}|\nabla u|^2\dd A
     +\int_{\Fill(E)}\frac{\dd A(\zeta)}{|\zeta-z|^2}.
\end{equation}
In particular $|\Phi'(z)|\le1$, with equality only if $F(w)=w+z$.
\end{lemma}
\begin{proof}
Set $v=\log|w|$, and integrate over $A_{\varepsilon,R}=\{\varepsilon<|w|<R\}\setminus J$, where the two boundary circles miss $J$. For these fixed radii the energy of $u+v$ is finite by change of variables under $F$:
solely a bounded region at positive distance from $z$ is integrated with weight $|\zeta-z|^{-2}$. The energy of $v$ on this annulus is $2\pi\log(R/\varepsilon)$.
The cross term is absolutely integrable by Cauchy--Schwarz. Almost every radial ray misses $J$, and polar integration gives
\[
 \int_{A_{\varepsilon,R}}\nabla u\cdot\nabla v\dd A
 =\int_0^{2\pi}\bigl(u(Re^{it})-u(\varepsilon e^{it})\bigr)\dd t
 \longrightarrow -2\pi\log|F'(0)|.
\]
The images of the two circles are Jordan curves. The outer curve lies in a circular band centered at $z$ of radii $R+O(1)$. The inner curve lies in a circular band of radii $\varepsilon |F'(0)|(1+O(\varepsilon))$. Both curves enclose $z$. For small $\varepsilon$ and large $R$, the region between their images contains the filled compact complement. The image of $A_{\varepsilon,R}$ is that region with $E^\sharp$ removed. The two circular bands have weighted area $o(1)$ in the joint limit. Thus
\[
 \int_{A_{\varepsilon,R}}|\nabla(u+v)|^2\dd A
 =2\pi\log(R/\varepsilon)-2\pi\log|F'(0)|
     -\int_{\Fill(E)}|\zeta-z|^{-2}\dd A(\zeta)+o(1).
\]
Singleton components have zero area, so the expression uses the original filled compact set. Expanding the square gives
\eqref{eq:energy} by exhaustion. If $\log|F'(0)|=0$, then $u$ is constant and tends to zero at infinity. The holomorphic function $(F(w)-z)/w$ has constant modulus $1$, so its normalization at infinity gives $F(w)=w+z$.
\end{proof}

\subsection{The two traces on a slit}

Let $J$ be such a union, with nondegenerate components $z_j(s)=\tau_j s$, $s\in I_j=[a_j,b_j]$, where $0<a_j<b_j$ and $|\tau_j|=1$. Degenerate slits are removable
and may be deleted. Set $p_j^2=\tau_j$, $|p_j|=1$, and write $z(x)$ and $p(x)$ for these functions on the disjoint union of the intervals, with measure $\dd\nu(s)=\dd s/\pi$. Use the complex Hilbert space $\mathcal H=L^2(\dd\nu)$. Let $C$ denote complex conjugation on $\mathcal H$ and write $L^T=CL^*C$ for a complex-linear operator $L$. Define the complex-linear operators
\begin{equation}\label{eq9}
 (Bu)(x)=\operatorname{pv}\int
   \frac{p(x)p(y)}{z(y)-z(x)}u(y)\dd\nu(y),\;K=-iB,\; M=iB .
\end{equation}
The operator $B$ is obtained by multiplying on both sides of the Cauchy--Hilbert transform on the arcs by $p$.
For the unweighted transform on a finite union of smooth arcs, see Murai \cite[Section~2, p.~315]{M1990}.
The diagonal blocks of $B$ are the finite Hilbert transforms with kernel $(y-x)^{-1}$ and measure $\dd\nu(y)$. Their boundedness follows by extension by zero and the $L^2$ bound for the Hilbert transform on the line \cite[Section~5.1]{Gr2014}. Each cross block is smooth on a neighborhood of its parameter rectangle, since the closed slits are disjoint. Thus $B$ is bounded and $B^T=-B$.

For the traces of $q_J$ and $h_J$, we use the Koebe map $L$ from the slit exterior
to a circular exterior with $L(w)=w+O(1)$ at infinity \cite[Theorem~2.4]{Y2018}.
To obtain the square root of $L'$ normalized at infinity, choose a smooth positively oriented Jordan loop around each slit. The side of the loop leading to infinity maps to the side of the image loop leading to infinity, since $L(\infty)=\infty$.
Conformality then preserves the orientation of the loop, so the tangent turning numbers of the loop and its image agree. Their difference is the winding number of $L'$ along the loop.
These loops generate the periods, so $L'$ has a holomorphic logarithm and hence a square root normalized at infinity. Pulling back the circular extremizers by multiplication by this square root gives $q_J,h_J$. Indeed, Lemma~\ref{lem2} and conformal invariance of the normalized
Ahlfors function identify the modulus of the pulled-back $q$ at every exterior point by Lemma~\ref{lem1}. Its normalization at infinity then identifies the function. The square-root pullback
preserves the boundary identity on each face. At each tip $a$, a local branch of $(w-a)^{1/2}$ opens a slit neighborhood into a half-disk, up to rotation.
By the prime-end boundary correspondence, $L$ extends continuously to the straightened diameter, and its values there lie on the boundary circle corresponding to the slit.
Schwarz reflection across the diameter and that circle \cite[p.~186]{N1952} extends $L$ holomorphically in the square-root coordinate. Consequently $L'(w)=O(|w-a|^{-1/2})$ and the pulled-back extremizers are $O(|w-a|^{-1/4})$. They extend analytically along each open face. Their traces are in $L^2(\dd\nu)$, and the integrals over the small endpoint circles vanish when the Cauchy contours collapse. The Sokhotski--Plemelj formula \cite[Section~17]{Mu1953} identifies the mean of the two traces with the principal value of the Cauchy integral.
For $L^2$ densities on a slit one can complete the slit to a smooth Jordan curve, extend the density by zero, and apply \cite[Theorem~2.3]{Le2018}.

For the side on the left of the directed ray $\tau_j$ write $Q_+=p_jq_+$ and $H_+=p_jh_+$. Use $Q_-,H_-$ on the right. The unit tangents in the compact-boundary orientation on the two faces are $-\tau_j$ and $\tau_j$, respectively.
The boundary relation in Lemma~\ref{lem1}, written with the exterior on the left by Murai \cite[p.~153, equation~(3)]{M1992}, gives in our orientation
$H_+=-iC Q_+$ and $H_-=iC Q_-$.
Let $u=(Q_++Q_-)/2$ and $m=(Q_+-Q_-)/2$. The mean and half-jump of the $H$ traces are $-iCm$ and $-iCu$. In particular,
\[
 q_+-q_-=\frac{2m}{p},\;\;\;h_+-h_-=-\frac{2i\,Cu}{p}.
\]
Define $P:\mathbb C\to\mathcal H$ by $P\omega=p\omega$ and write $P_0=P1$.
The exterior Cauchy formula for $f=q_J,h_J$ has a minus sign when the contours are oriented as boundaries of the compact sets.
Collapsing the contours to the two faces therefore gives
\[
 f(w)-f(\infty)
   =\frac{1}{2\pi i}\sum_j\int_{I_j}
     \frac{(f_+(s)-f_-(s))\tau_j}{z_j(s)-w}\dd s .
\]
Substituting the two jumps, using $\tau_j=p_j^2$ and $\dd\nu=\dd s/\pi$, we obtain
\begin{equation} \label{eq10}
 q_J(w)=1-i\int\frac{p(y)m(y)}{z(y)-w}\dd\nu(y),\;\;\;h_J(w)=-\int\frac{p(y)(Cu)(y)}{z(y)-w}\dd\nu(y).
\end{equation}
Murai \cite[p.~319, equation~(10) and the preceding display]{M1990} gives the corresponding Cauchy representations on smooth arcs.
In his formulas the density corresponding to our variables is $-ipm$, the normalized Garabedian square root is $q_J$, and the numerator in the formula for the Ahlfors function is $h_J$.
Taking the means of their two traces gives $u=P_0+Km$ and $Cm=KCu$. In particular $m=\overline K u=-K^*u$. Here $C$ is the pointwise complex conjugation operator.

\subsection{Identification of a relative determinant}

View $\mathcal H$ also as a real Hilbert space, with the real part of its complex inner product, and put $D=MC=-KC$. Since $M^T=-M$, one has $\overline M=-M^*$ and
$D^*=(CM^*C)C=-D$ for the real adjoint. Also $(I+D)(I-D)=(I-D)(I+D)=I+MM^*$, so $I+D$ is invertible. The trace relations give $u+Mm=P_0$ and $Cm+MCu=0$. It follows that, for $\omega\in\mathbb C$,
\begin{equation}\label{eq11}
   (I+D)^{-1}P\omega=u\omega+(Cm)\overline\omega .
\end{equation}

Let $W$ be multiplication by $z(x)$, and set $D'=W D\overline W^{-1}$. The multiplier $W$ is bounded and invertible since the slits are bounded away from zero. Let $G:\mathbb C\to\mathcal H$ be given by
$G\omega=(p/z)\omega$, and let $J_0\omega=i\overline\omega$. We use complex adjoints for complex-linear maps, which agree with their real adjoints. Multiplication of the Cauchy kernel
by $z(x)/z(y)$ gives
\begin{equation}\label{eq12}
 W B W^{-1}=B-P\ell,\;
 \ell(v)=\int\frac{p(y)v(y)}{z(y)}\dd\nu(y),\;
 D'=D-PJ_0G^* .
\end{equation}
Here $G^*v=\overline{\ell(Cv)}$. Thus \eqref{eq11} gives
\[
 (I-J_0G^*(I+D)^{-1}P)\omega
   = (1-i\ell(m))\omega-i\ell(Cu)\overline\omega
   = q_J(0)\omega+i h_J(0)\overline\omega .
\]
The real determinant of $\omega\mapsto a\omega+b\overline\omega$ is $|a|^2-|b|^2$. For the map in the preceding display, Smith's slope formula (Lemma~\ref{lem1}) identifies this determinant with $a_J(0)$. To apply Sylvester's determinant identity, regard $\C$ as a real plane and put
$R_0=J_0G^*(I+D)^{-1}:\mathcal H\to\C$, with both spaces viewed over $\R$. Then $P:\C\to\mathcal H$ is the map in the opposite direction, and \eqref{eq12} gives
$(I+D')(I+D)^{-1}=I-PR_0$. The range of $PR_0$ is contained in the finite-dimensional space $S=\operatorname{ran}P$,
so the Fredholm determinant of $I-PR_0$ is the determinant of its restriction to $S$ \cite[Chapter~3]{Si2005}. Applying the ordinary Sylvester identity to
$P:\C\to S$ and $R_0|_S:S\to\C$ then gives
\[
 \det_{\R}(I-PR_0)=\det_{\R}(I-R_0P).
\]
The determinant on the right is the real-plane determinant evaluated above. Thus
\begin{equation}\label{eq13}
 a_J(0)=\det_{\mathbb R}\big((I+D')(I+D)^{-1}\big).
\end{equation}
Here we use the trace ideals and the Fredholm determinant in \cite[Chapters~2--3]{Si2005}. The ideals $\mathfrak S_1$ and $\mathfrak S_2$ consist of compact operators with summable and square-summable singular values, respectively. For $B\in\mathfrak S_1$, $\det(I+B)$ is the product of $1+\lambda_k(B)$ with eigenvalues counted with algebraic multiplicity. On a real Hilbert space we take this determinant
on the complexification. We use multiplicativity on $I+\mathfrak S_1$, invariance under bounded similarity, and continuity in trace norm.

\subsection{Passage to an orthogonal projection}

Write $W=HU$, where $H$ is multiplication by $|z|$ and $U$ is multiplication by $\tau_j$ on $I_j$. Then
\[
 D'=H\widetilde D H^{-1},\;
 \widetilde D=U D U=N C,\; N=U M U^{-1}.
\]
We use the following observation about integral kernels. Let $F(s,t)$ be smooth on a neighborhood of a compact parameter rectangle $I_j\times I_k$. Choose a smooth cutoff supported in that neighborhood and equal to one on the parameter rectangle. Extend the product by zero and periodize it on a larger rectangle with periods $L_s,L_t$. Write $\widetilde F$ for this periodic extension and $c_{mn}$ for its Fourier coefficients, so that on $I_j\times I_k$,
\[
 F(s,t)=\sum_{m,n\in\mathbb Z}c_{mn}e^{2\pi i ms/L_s}e^{2\pi i nt/L_t}.
\]
In fact, integrating by parts on the two circles as in \cite[Section~3.3]{Gr2014} gives
\[
 |c_{mn}|\le \frac{C_F}{(1+m^2)(1+n^2)}.
\]
Here $C_F$ is the mean over a period rectangle of the absolute value of
\[
 \bigl(1-(L_s/2\pi)^2\partial_s^2\bigr)
 \bigl(1-(L_t/2\pi)^2\partial_t^2\bigr)\widetilde F .
\]
The series therefore converges absolutely and uniformly, and its sum agrees with $\widetilde F$ by uniqueness of Fourier coefficients.
On $L^2(I_k,\dd\nu)$ the $(m,n)$ term of the integral operator is
\[
 f\longmapsto c_{mn}e^{2\pi i ms/L_s}
                 \int_{I_k}e^{2\pi i nt/L_t}f(t)\dd\nu(t).
\]
This operator has rank at most one and trace norm $|c_{mn}|\sqrt{\nu(I_j)\nu(I_k)}$, by the formula for the norm of a rank-one operator.
Consequently the operator series converges in trace norm, and uniform convergence of the kernels identifies its limit with the integral operator for $F$.
If either parameter already runs on a circle, we use its given period and apply the cutoff and extension only in the interval variables.

For $j\ne k$, the block of $N-M$ from $I_k$ to $I_j$ has kernel
\[
 i\left(\frac{\tau_j}{\tau_k}-1\right)\frac{p_jp_k}{\tau_k t-\tau_j s}.
\]
The denominator is bounded away from zero because the closed slits are disjoint. The kernel is thus smooth on a neighborhood of $I_j\times I_k$.
The diagonal blocks of $N-M$ vanish since $U$ is constant on each $I_j$. There are only finitely many blocks, so $N-M\in\mathfrak S_1$.

Since $H$ is real multiplication, $[H,D]=[H,M]C$. The block of $[H,M]$ from $I_k$ to $I_j$ has kernel
\[
 i p_jp_k\frac{s-t}{\tau_k t-\tau_j s}.
\]
For $j=k$ this is the constant $-i$ after cancellation. The truncated kernels converge to this constant in $L^2$ of the product measure, so the resulting operator has rank at most one.
The other blocks are smooth by the same separation argument.
Thus $[H,M]\in\mathfrak S_1$. If $V=[H,M]$, the complexification of the real operator $VC$ is
$\begin{pmatrix}0&V\\\overline V&0\end{pmatrix}$ on $\mathcal H\oplus\mathcal H$, where $\overline V=CVC$ has the same singular values as $V$.
The complexified operator is trace class, as required for the real determinant.

For $A_s=I+sD$, $0\le s\le1$, write $C_s=H A_s H^{-1}A_s^{-1}$ and $B_s=A_s^*$. The identity $C_s-I=s[H,D]H^{-1}A_s^{-1}$ shows that $C_s-I$ is trace class on the complexification and varies continuously in trace norm, since $A_s^{-1}$ varies continuously in operator norm. Conjugation
by multiplication by $i$ commutes with $H$ and exchanges $A_s$ and $B_s$, so $C_s$ and $H B_s H^{-1}B_s^{-1}$ have the same determinant. On the other hand,
$C_s^*=B_s^{-1}H^{-1}B_s H
  =(B_s H)^{-1}(H B_s H^{-1}B_s^{-1})^{-1}(B_s H).$
The determinant on the real space is real and equals that of the adjoint. By similarity invariance and multiplicativity \cite[Chapter~3]{Si2005},
$\det_{\R} C_s=(\det_{\R} C_s)^{-1}$. Continuity in trace norm and $C_0=I$ give $\det_{\R} C_s=1$.
In \eqref{eq13}, factor the relative operator as $ H\big((I+\widetilde D)(I+D)^{-1}\big)H^{-1} C_1 .$
This gives
\begin{equation}\label{eq14}
 a_J(0)=\det_{\mathbb R}\big((I+\widetilde D)(I+D)^{-1}\big).
\end{equation}

Complexification sends $I+MC$ to $\begin{pmatrix}I&M\\\overline M&I\end{pmatrix}$.
Factor it as
\[
 \begin{pmatrix}I&M\\0&I\end{pmatrix}
 \begin{pmatrix}I-M\overline M&0\\0&I\end{pmatrix}
 \begin{pmatrix}I&0\\\overline M&I\end{pmatrix},
\]
and use the identical factorization for $N$. The relative upper and lower triangular factors have determinant one because $N-M$ is trace class and their off-diagonal
perturbations square to zero. Denote the three factors by $L_M,S_M,R_M$ and use the same notation with subscript $N$ for the second product. With $\operatorname{Ad}_P(V)=PVP^{-1}$,
the relative product is
\[
 (L_NS_NR_N)(L_MS_MR_M)^{-1}
   =(L_NL_M^{-1})\,
     \operatorname{Ad}_{L_M}(S_NS_M^{-1})\,
     \operatorname{Ad}_{L_MS_M}(R_NR_M^{-1}).
\]
The middle factor differs from the identity by a trace-class operator since $N-M$ is trace class. Using $\overline M=-M^*$ and setting $A=I+MM^*$, we obtain
\begin{equation}\label{eq15}
    a_J(0)=\det_{\mathbb C}\big((I-N\overline N)A^{-1}\big).
\end{equation}

Put $X=N^*$, $Y=-\overline N$, $A_X=I+X^*X$, and $A_Y=I+Y^*Y$. Then
$I+X^*Y=I-N\overline N,\;A_X=U A U^{-1},\;A_Y=U^{-1} A U .$ The relative normalization determinants cancel: $A_YA^{-1}=U^{-1}(A_XA^{-1})^{-1}U$. Both are identity plus trace class because $N-M$ is trace class. Hence their determinants have product one.
The isometries
\[
 Q_X=\binom{I}{X}A_X^{-1/2},\qquad
 Q_Y=\binom{I}{Y}A_Y^{-1/2}
\]
have ranges the graphs of $X$ and $Y$, respectively. Thus $\Pi_X=Q_XQ_X^*$ is the orthogonal projection onto the graph of $X$
\cite[Theorem~4]{Sto1951}. The positive contraction $T=Q_Y^*\Pi_XQ_Y$ is similar to $B_1^*A_X^{-1}B_1 A_Y^{-1}$, where $B_1=I+X^*Y$. The latter product is
$(B_1^* A^{-1})(A A_X^{-1})(B_1 A^{-1})(A A_Y^{-1}).$
Each factor differs from $I$ by a trace-class operator. The first and third have conjugate determinants by adjoint and similarity, and the other two have determinants with product
$1$. By multiplicativity \cite[Theorem~3.5]{Si2005}, \eqref{eq15} gives
\begin{equation}\label{eq16}
        a_J(0)^2=\det_{\mathbb C}(Q_Y^*\Pi_XQ_Y)\le1 .
\end{equation}
Since $T$ is similar to the product above, $T-I$ is trace class. The eigenvalues of $T-I$ lie in $[-1,0]$ because $0\le T\le I$.
Each factor in the spectral product for $\det_{\C}T$ therefore belongs to $[0,1]$, which proves the inequality \cite[Chapter~3]{Si2005}. Monotonicity of capacity
gives $a_J(0)\ge0$, hence $a_J(0)\le1$.

\subsection{Analytic curves and the general slope bound}

\begin{proposition}\label{prop1}
Let $K$ be a finite union of mutually disjoint analytic Jordan compacts and regular compact analytic arcs. Here a Jordan compact is a Jordan curve together with its bounded interior.
For every finite $z\in\Omega_K$, the graph determinant in \eqref{eq16} gives $a_K(z)^2$ on the boundaries and arcs of $K$. Translate $z$ to zero, parametrize the integration curves by
arclength $s$, and in \eqref{eq9} use $p(s)^2=\zeta'(s)$ and $U(s)=\zeta(s)/|\zeta(s)|$, where $\zeta$ is the resulting boundary or arc coordinate. Equality $a_K(z)=1$ for nonempty $K$ requires that all the integration curves lie on a single straight line through $z$, and this condition suffices.
\end{proposition}
\begin{proof}
Orient each closed curve counterclockwise and set the left traces of $q_K,h_K$ equal to zero there, using the exterior values on the right. The corresponding trace identity is
$H_-=i\overline{Q_-}$, while the zero traces satisfy the identity on the left. At an arc the two traces satisfy the relations used for slits. In fact, an analytic coordinate
straightens the arc near each endpoint, after which the square-root chart and circular uniformization give the same $O(\dist(\zeta,\text{tip})^{-1/4})$ estimate. The Cauchy formula with all these traces gives \eqref{eq10} and the same mean and jump equations. One may choose a measurable square root $p$ on a closed curve, since all occurrences of $p$ in the kernel are separate bounded multipliers in the two variables.

On a single parameter interval, write $v=t-s$. On a closed curve, use the same notation in a local lift of the parameter circle near its diagonal, with parameters taken modulo the length $l$.
In either case, $\zeta(s+v)-\zeta(s)=v A(s,v),\;\;A(s,v)=\int_0^1\zeta'(s+\lambda v)\dd\lambda .$
The function $A$ is smooth and $A(s,0)=\zeta'(s)\ne0$. Thus $1/A(s,v)-1/\zeta'(s)$ is divisible by $v$ with smooth quotient.
For an arc this shows that the unweighted Cauchy kernel differs from $1/(\zeta'(s)(t-s))$ by a smooth kernel.
For a closed curve the model kernel is
$\zeta'(s)^{-1}(\pi/l)\cot(\pi(t-s)/l)$, since $(\pi/l)\cot(\pi v/l)-1/v$ extends smoothly across $v=0$.
The remainder extends smoothly across the parameter diagonal on the torus. Away from the diagonal the kernels are smooth by injectivity of each parametrization.
On blocks between distinct curves the Cauchy kernel is smooth by separation. Using symmetric arclength truncations on each diagonal block, the real-line and periodic Hilbert transforms
\cite[Sections~5.1 and~4.1.2]{Gr2014} give the $L^2$ limits and boundedness, also after multiplication by $p$ on both sides.
Truncation at chord length $\varepsilon$ gives the same operator in the $L^2$ limit. Indeed, unit speed gives
$|\zeta(s+v)-\zeta(s)|=|v|+O(v^2)$ uniformly near the diagonal. By compactness and injectivity, for small $\varepsilon$ the two truncations can differ only where $|v|=\varepsilon+O(\varepsilon^2)$.
There the kernel is $O(\varepsilon^{-1})$, and the measure of the strip in each variable with the other fixed is $O(\varepsilon^2)$. Schur's test \cite[Appendix~A]{GrM2014} bounds the difference of the operators by $O(\varepsilon)$.
The transpose identity follows by passing to the limit in the antisymmetry of the symmetrically truncated kernels.

Put $\rho(s)=|\zeta(s)|$ and $\eta(s)=\zeta(s)/\rho(s)$, so that $H$ and $U$ are multiplication by $\rho$ and $\eta$, respectively.
The marked point lies off the curves, so $\rho$ is positive and both functions are smooth on the parameter intervals or circles.
After removing the factors $ip(s)p(t)$, the kernels of $N-M$ and $[H,M]$ are, respectively,
\[
 \frac{\eta(s)/\eta(t)-1}{\zeta(t)-\zeta(s)}
 \quad\text{and}\quad
 \frac{\rho(s)-\rho(t)}{\zeta(t)-\zeta(s)} .
\]
On a block belonging to a single curve their smooth extensions across the diagonal follow from the formulas
\[
 \frac{\eta(s)/\eta(s+v)-1}{\zeta(s+v)-\zeta(s)}
   =-\frac{\int_0^1\eta'(s+\lambda v)\dd\lambda}{\eta(s+v)A(s,v)},\;\;
 \frac{\rho(s)-\rho(s+v)}{\zeta(s+v)-\zeta(s)}
   =-\frac{\int_0^1\rho'(s+\lambda v)\dd\lambda}{A(s,v)} .
\]
On blocks between distinct curves the displayed quotients are smooth on the entire product of the parameters. The cancelled kernels are bounded.
Hence their truncations converge to them in $L^2$ of the product measure, and the corresponding operators converge in Hilbert--Schmidt norm to the ordinary integral operators.
The Fourier argument above gives trace-class operators on all these blocks. The bounded multiplication operators by $p$ on the source and target preserve trace class.
As in the slit case, $[H,D]=[H,M]C$ is trace class after complexification. All steps from \eqref{eq11} through \eqref{eq16} now apply, including the cancellation of the two normalization determinants.

If $a_K(z)=1$, then the positive contraction $T=Q_Y^*\Pi_XQ_Y$ has determinant one. Since $T-I$ is trace class, the spectral product gives $T=I$.
The identity
\[
 \langle (I-T)v,v\rangle=\|(I-\Pi_X)Q_Yv\|^2
\]
shows that the graph of $Y$ lies in the graph of $X$, so $Y=X$. Equivalently, $U^2$ commutes with $M^*$.
Restrict the input and output of this commutator to two parameter neighborhoods with disjoint closed images. On their product it is an ordinary Hilbert--Schmidt integral operator with kernel
$(U(s)^2-U(t)^2)M^*(s,t)$, where $M^*(s,t)$ denotes the kernel of $M^*$ away from the diagonal.
A zero Hilbert--Schmidt integral operator has kernel zero almost everywhere. A countable collection of such products covers all off-diagonal pairs, and $M^*(s,t)$ is nonzero there.
Consequently $U(s)^2=U(t)^2$ for almost every pair of parameters, and continuity gives the equality everywhere. This is collinearity through the marked point. The converse follows from the same commutation and the determinant formula.
\end{proof}

\begin{proof}[Proof of Theorem~\ref{thm1}] 
For an analytic exterior, Proposition~\ref{prop1} and the projection inequality in \eqref{eq16} give $a_E(z)\le1$. For general $E$, take decreasing analytic outer approximants to
$\Fill(E)$. Their slopes converge at every $z\in\Omega_E$ by Smith's convergence theorem \cite[Theorem~3.1]{S1974}, as in Lemma~\ref{lem1}. The lemma also gives strict positivity, and a singleton attains the bound.

For later use, in finite connectivity we apply the determinant formula to the radial-slit exterior of Section~\ref{sec3}. As $a_E(z)=a_{E^\sharp}(z)$, Lemmas \ref{lem2} and \ref{lem3} give the exact factorization
\begin{equation}\label{eq:full-factorization}
 a_E(z)=
 \exp\!\left(-\frac1{2\pi}
   \left[\int_{\C\setminus J}|\nabla u|^2\dd A
        +\int_{\Fill(E)}\frac{\dd A(\zeta)}{|\zeta-z|^2}\right]\right)
 \det\nolimits_{\C}(Q_Y^*\Pi_XQ_Y)^{1/2}.
\end{equation}
For $J=\varnothing$ the determinant is $1$.
\end{proof}

\begin{corollary}\label{cor1}
Suppose that $\Omega_E$ is finitely connected and that $z\in\Omega_E$ is finite. Then $a_E(z)=1$ if and only if all the non-singleton components of $\Fill(E)$ are straight segments on a single line through $z$. If there are no such components, the condition is vacuous.
\end{corollary}
\begin{proof}
The slope is unchanged on passing to $\Omega_E^\sharp$. By \eqref{eq:full-factorization}, equality requires $F(w)=w+z$ on $\Chat\setminus J$ and equality in \eqref{eq16}.
A positive contraction differing from the identity by a trace-class operator has determinant one only if it is the identity. In \eqref{eq16} this says that the graph of $Y$ is contained in the graph
of $X$, and hence $Y=X$. Equivalently $U^2$ commutes with $M^*$. Every cross block of $M^*$ between two nondegenerate slits is nonzero, so this happens exactly when $\tau_j^2$ is independent of $j$. Thus the slits are collinear. Conversely, for collinear segments through the marked line, translation already puts the exterior in radial-slit form, the two graphs coincide, and \eqref{eq:full-factorization} gives equality.
\end{proof}

If an arbitrary compact set is contained in a line through $z\in\Omega_E$, approximate it from above in that line by finite unions of closed intervals missing $z$. The corollary and the slope convergence after Lemma~\ref{lem1} give $a_E(z)=1$ also in this case, in agreement with Smith \cite[Corollary~16.3]{ST1974}.

\begin{remark}\label{rem1}
The graph formula also bounds the deficit away from the collinear configuration. Define
\[
 V=(I+XX^*)^{-1/2}(Y-X)A_Y^{-1/2}.
\]
The isometry $\binom{-X^*}{I}(I+XX^*)^{-1/2}$ onto the orthogonal complement of the graph of $X$ gives $Q_Y^*(I-\Pi_X)Q_Y=V^*V$. By \eqref{eq16} and $1-t\le e^{-t}$ for $0\le t\le1$,
\[
 a_J(0)\le \exp\!\left(-\frac{\|Y-X\|_{\mathfrak S_2}^2}
                                  {2(1+\|M\|^2)^2}\right).
\]
The diagonal blocks of $Y-X$ vanish and the modulus of its $(i,j)$ cross kernel is $|\tau_i^2-\tau_j^2|/|\tau_i s-\tau_j t|$.
\end{remark}

\section{Comparison of finite increments}\label{sec4}
In this section, we present the complete proof of Theorem \ref{thm2}. To do so, we first establish the following two results.
\begin{proposition}\label{prop3}
Assume \hyperref[eqS]{(S)}. If $C\subset B$ are compact and $C$ is a continuum, then $a_B(z)\le a_C(z)$ at every finite point of
$\Omega_B$.
\end{proposition}
\begin{proof}
If $C$ is a singleton, then $a_C(z)=1$ by removability, and \hyperref[eqS]{(S)} gives $a_B(z)\le a_C(z)$.
Otherwise set $r=\gamma(C)>0$ and let $U$ map $\Omega_C$ onto the exterior of the closed disk of radius $r$, with $U(w)=w+O(1)$ at infinity.
The explicit disk extremizers $q=1$ and $h=r/w$, and Lemma~\ref{lem2}, give $a_C(z)=|U'(z)|(1-r^2/|U(z)|^2)$. Put $\theta=\arg U(z)$.
The Joukowski map \cite[Example~3.6, p.~185]{Sa2010}  $\Phi(w)=U(w)+r^2e^{2i\theta}/U(w)$ is univalent on $\Omega_C$: the difference of two images under
$u\mapsto u+r^2e^{2i\theta}/u$ is $(u-v)(1-r^2e^{2i\theta}/(uv))$, and the second factor cannot vanish when $|u|,|v|>r$.
Direct differentiation gives $|\Phi'(z)|=a_C(z)$. Restrict $\Phi$ to $\Omega_B$ and denote the compact complement of its image by $B'$.
By Lemma~\ref{lem2} and \hyperref[eqS]{(S)}, $a_B(z)=|\Phi'(z)|a_{B'}(\Phi(z))\le|\Phi'(z)|=a_C(z)$.
\end{proof}

\begin{lemma}\label{lem4}
Suppose that one analytic Jordan boundary of an analytic exterior moves smoothly with exterior normal speed $v$. Let $n$ be the unit normal pointing into the exterior. If $K_t$ is the resulting compact set,
then at the stated time, with the integral taken over the moving boundary,
\begin{equation}\label{eq8}
  \frac{\dd}{\dd t}\gamma(K_t)
     =\frac1{4\pi}\int v\,\partial_n a_{K_t}\dd s .
\end{equation}
\end{lemma}
\begin{proof}
Fix the stated time and use its extremizers $q,h$ in the two problems \eqref{eq5} at a nearby time. Analytic extension makes
them admissible on both sides of the fixed time. If $k$ is the derivative of the argument of the positively oriented tangent with respect to arclength, the derivatives of the upper and lower trial values are
\[
  \frac1{2\pi}\int v(\partial_n|q|^2+k|q|^2)\dd s,\;\;\;-\frac1{2\pi}\int v(\partial_n|h|^2+k|h|^2)\dd s .
\]
With the compact-boundary orientation we have $n=-iT$ and $k=\frac{\dd}{\dd s}\arg T$. On the boundary, $qhT=i|q|^2$ and $|h|=|q|$.
At points where $qh\ne0$, the Cauchy--Riemann equations give
\[
 \partial_n\log|qh|=\frac{\dd}{\dd s}\arg(qh)
   =-\frac{\dd}{\dd s}\arg T=-k .
\]
Thus $\partial_n(|q|^2+|h|^2)=-2k|q|^2$, with the identity extending by continuity at zeros. Both derivatives above equal the right-hand side $d_0$ of \eqref{eq8} at the fixed time $t_0$.
For either sign of $\varepsilon$, the lower and upper trial values at $t_0+\varepsilon$ bound $\gamma(K_{t_0+\varepsilon})$ and both equal $\gamma(K_{t_0})+\varepsilon d_0+o(|\varepsilon|)$.
This gives the two-sided derivative in \eqref{eq8}.
\end{proof}

\subsection{Separated continua}

Suppose first that $C_0,C_1$ are strictly nested analytic Jordan compacts and that the filled analytic compact $E$ is disjoint from $C_1$. A conformal map of the annular region between their
boundaries onto a round annulus gives a family $C_t$, $0\le t\le1$, expanding from $C_0$ to $C_1$ with nonnegative normal speed. The boundaries remain analytic, and their analytic parametrizations extend
locally past both end times by reflection. On the exterior side of $\partial C_t$, Proposition~\ref{prop3} yields $a_{E\cup C_t}\le a_{C_t}$. Both slopes tend to zero on the boundary and extend smoothly to it, by Lemma~\ref{lem1}. Therefore their outward normal derivatives satisfy the same inequality. Lemma~\ref{lem4} gives
\[
   \frac{\dd}{\dd t}\gamma(E\cup C_t)
   \le \frac{\dd}{\dd t}\gamma(C_t) .
\]
Thus $g(t)=\gamma(E\cup C_t)-\gamma(C_t)$ satisfies $g'(t)\le0$. The mean value theorem gives $g(1)\le g(0)$, which is \eqref{eq3}.

For disjoint $E,C_1$, retain only $E'=E\cap\Omega_{C_1}$. This is compact by positive separation and $\partial\Omega_{C_1}\subset C_1$. The terminal union has the same capacity after this replacement, while the initial union can have smaller capacity. Thus it suffices to prove the upper bound with $E'$. If $E'=\varnothing$, the comparison with $E'$ is an equality.
Otherwise, fill $C_0,C_1$. The connected filled $C_1$ is contained in one component of $\Chat\setminus E'$.
If this component is bounded, adjoining either continuum leaves $\gamma(E')$ unchanged. Otherwise $\Fill(E')$ is disjoint from the filled $C_1$. Outer Riemann-map levels
approximate each filled continuum by decreasing analytic Jordan compacts, using disks for a singleton. For each approximant to $C_1$, choose the approximant to $C_0$ inside its interior. The analytic neighborhoods of $\Fill(E')$ and $C_1$ can be kept disjoint by positive separation. Now use continuity from above in the analytic inequality.

\subsection{Intersections}

We use the area Cauchy transform as in \cite{V1967,Ga1972}. For any bounded compactly supported density $\beta$ write $\Fcal_{A}\beta(z)=\pi^{-1}\int\beta(w)/(z-w)\dd A(w)$. It solves
$\dbar\Fcal_A\beta=\beta$ in distributions. If a summand of $\beta$ is supported in a disk of radius $s$ with modulus bounded by $M$, its Cauchy correction has supremum norm at most $CMs$
and coefficient of $1/z$ bounded by $CMs^2$. This follows by integrating $|z-w|^{-1}$ on the disk, and by expanding at infinity.

For the moment suppose that $E$ has finite analytic boundary, that $C_1$ is an analytic Jordan compact, and that their boundaries meet transversely at a finite set $P$. Put
\[
  K=E\cup C_1,\;\;E_\delta=\{w\in E:\dist(w,C_1)\ge\delta\},\;\;K_\delta=E_\delta\cup C_1.
\]
Any point of $\partial K$ exposed in $\C\setminus K_\delta$ lies on $\partial E\setminus C_1$ at distance less than $\delta$ from $C_1$. For every fixed neighborhood of $P$, all such points
lie in it when $\delta$ is sufficiently small. The assertion also holds when $P=\varnothing$.

For each $p\in P$, let $u_p$ map the exterior of $C_1$ to $\D$,  with $u_p(\infty)=0$ and $u_p(p)=1$. This map extends holomorphically near $p$. For $\alpha>0$ set $g_p=((1-u_p)/2)^\alpha$, using the branch in the right half-plane. Then $|g_p|\le1$, $g_p(\infty)=2^{-\alpha}$, and $g_p(w)=O(|w-p|^\alpha)$. Write $g=\prod_{p\in P}g_p$. Take the Ahlfors function $f$ of $K$ on its full complement, with value
zero on bounded components. Choose a smooth cutoff $\chi$ equal to zero on the $2\eta$-neighborhood of $P$ and equal to one outside its $3\eta$-neighborhood, with $|\nabla\chi|\le C/\eta$.
Fix $\alpha$ and $\eta$ and take $\delta$ so small that the newly exposed boundary lies in the $\eta$-neighborhood of $P$. There $\chi$ vanishes in a neighborhood, so $\chi fg$ extends by zero from
$\C\setminus K$ to $\C\setminus K_\delta$ and its distributional $\dbar$ derivative is just the ordinary cutoff term. Extend this derivative $\beta$ by zero to the whole plane.
The Cauchy correction $\Fcal_A\beta$ has supremum norm $O(\eta^\alpha)$ and $1/z$ coefficient $O(\eta^{1+\alpha})$. The corrected function $\chi fg-\Fcal_A\beta$ has vanishing distributional
$\dbar$ on $\C\setminus K_\delta$, so Weyl's lemma gives a holomorphic representative there with the same norm bound and Laurent coefficient. It yields
\[
 \liminf_{\delta\downarrow0}\gamma(K_\delta)
  \ge\frac{2^{-\alpha|P|}\gamma(K)-O(\eta^{1+\alpha})}
          {1+O(\eta^\alpha)} .
\]
Let $\eta\downarrow0$ and then $\alpha\downarrow0$. Monotonicity gives $\gamma(K_\delta)\to\gamma(K)$. This is a finite-interface version of the local gap-opening argument in
Smith \cite[Lemma~18.1]{ST1974}.

We can therefore apply the separated result with $E_\delta$ and the original nested continua. Use $\gamma(E_\delta\cup C_0)\le\gamma(E\cup C_0)$ and pass to the limit. For general $E,C_0,C_1$, fill the three sets and take decreasing analytic outer approximants $E_n$ to $\Fill(E)$. Use the outer levels of the modulus of the exterior Riemann map of the filled $C_1$ onto the exterior of the unit disk,
or the distance to $C_1$ when $C_1$ is a point. The limiting level is $1$ for the modulus and $0$ for the distance.
At stage $n$, choose a regular value on the portions of $\partial E_n$ in the exterior of $C_1$, with excess above the limiting level smaller than $1/n$ and smaller than at the previous stage.
Such a choice is possible because the critical values have measure zero.
The resulting analytic Jordan compacts $C_{1,n}$ decrease to the filled $C_1$, and their boundaries meet $\partial E_n$ transversely in finitely many points.
Choose $C_{0,n}$ by an outer Riemann-map level, or by a small disk for a singleton, close enough to the filled $C_0$ to lie inside $\operatorname{int}C_{1,n}$ and, for $n>1$, inside $\operatorname{int}C_{0,n-1}$.
The levels can be chosen to tend to the limiting level for $C_0$. Continuity from above proves \eqref{eq3}.

\begin{proof}[Proof of Theorem~\ref{thm2}]
The preceding argument proves \textup{(i)}$\Rightarrow$\textup{(iii)}. Taking $C_0$ to be a point of the disk $C_1$ gives \textup{(iii)}$\Rightarrow$\textup{(ii)}, since a point has no
effect on capacity. Finally $\gamma(\overline D(z,s))=s$. Divide the resulting increment bound by $s$ and let $s\downarrow0$ to get \textup{(ii)}$\Rightarrow$\textup{(i)}.
\end{proof}

\begin{corollary}\label{cor2}
Under \hyperref[eqS]{(S)}, for any compact $E$ and continuum $C$, $\gamma(E\cup C)\le\gamma(E)+\gamma(C)$. If $E\cap C$ is a continuum, then $\gamma(E\cup C)+\gamma(E\cap C)\le\gamma(E)+\gamma(C)$.
\end{corollary}
\begin{proof}
Use \eqref{eq3} with $C_1=C$ and take $C_0$ to be a singleton or $E\cap C$, respectively.
\end{proof}

\begin{corollary}\label{cor3}
Let $E\subset\C$ be compact and let $F$ be a compact subset of a straight line. Then
\begin{equation}\label{eq:collinear-increment}
 0\le\gamma(E\cup F)-\gamma(E)
 \le \tfrac14\Hm^1(F\setminus E)
 =\gamma(F)-\gamma(E\cap F).
\end{equation}
In particular, $\gamma(E\cup F)+\gamma(E\cap F)
\le\gamma(E)+\gamma(F)$.
\end{corollary}
\begin{proof}
The assertion is immediate if either set is empty. Let $\ell$ be a line containing $F$. Given $\varepsilon>0$, outer regularity of length on $\ell$ gives bounded relatively open intervals
$I_j\subset\ell$ covering $F\setminus E$ and satisfying
\[
 \sum_j |I_j|\le \Hm^1(F\setminus E)+\varepsilon .
\]
For $\delta>0$, let $E_\delta$ be the closed $\delta$-neighborhood of $E$. The intervals $I_j$ and the relatively open set $\ell\cap\operatorname{int}E_\delta$ cover $F$. By compactness,
a finite set of intervals indexed by $J$ suffices together with $\ell\cap\operatorname{int}E_\delta$. Adjoin the closed intervals $\overline I_j$, $j\in J$, to $E_\delta$ one at a time.
Corollary~\ref{cor2}, with the slope bound from Theorem~\ref{thm1}, bounds each increment by $\gamma(\overline I_j)=|I_j|/4$. Monotonicity therefore gives
\[
 \gamma(E\cup F)\le\gamma(E_\delta)+\tfrac14\sum_{j\in J}|I_j|
 \le\gamma(E_\delta)+\tfrac14\bigl(\Hm^1(F\setminus E)+\varepsilon\bigr).
\]
Let $\delta\downarrow0$ using continuity from above and then let $\varepsilon\downarrow0$. The lower bound follows from monotonicity. Since $F$ and $E\cap F$ are compact subsets of $\ell$,
Pommerenke's formula \cite{P1960} gives the final equality in \eqref{eq:collinear-increment}.
\end{proof}

\section{Length, Cauchy capacity, and null sets}\label{sec5}

For $A\subset\C$, use the disk-cover normalization of one-dimensional Hausdorff content with radii as
covering costs \cite[Section~2]{Par2022}, and define the disk content by
\[
  \mathcal D^1(A)=\inf\left\{\sum_j r_j:
                      A\subset\bigcup_j D(z_j,r_j)\right\},
\]
with countable covers allowed. In particular $\mathcal D^1(A)\le\Hm^1(A)$, since a nonempty set of diameter $d$ can be covered by an open disk of radius arbitrarily close to $d$.

\begin{proposition}\label{prop4}
For compact $E,F\subset\C$,
\begin{equation}\label{eq17}
 \gamma(E\cup F)-\gamma(E)\le\mathcal D^1(F\setminus E).
\end{equation}
\end{proposition}
\begin{proof}
Fix an open disk cover $D(z_j,r_j)$ of $F\setminus E$. Write $E_\varepsilon$ for the closed $\varepsilon$-neighborhood of $E$. The covering disks and $\operatorname{int} E_\varepsilon$ have
a finite subcover of $E\cup F$. Adjoin the selected closed disks one at a time to $E_\varepsilon$. By Theorem~\ref{thm2}, the increment is at most $\sum_j r_j$. Now let $\varepsilon\downarrow0$ and take the infimum.
\end{proof}

\begin{proof}[Proof of Theorem~\ref{thm3}]
Write $L=\Hm^1(R\setminus K_0)$ and fix $\varepsilon>0$. Assign points of $R\setminus K_0$ to the first one of the
countable embedded arcs containing them, discarding a null set. Take finitely many of these disjoint Borel pieces leaving mass less than $\varepsilon$, and compact subsets of the chosen pieces 
leaving additional mass less than $\varepsilon$. By the area formula, the injective unit-speed parametrization of each carrier preserves $\Hm^1$ on Borel subsets of its parameter interval.
In the arclength coordinate, cover the selected compact subset by finitely many intervals whose total length is at most the measure of the piece plus its allotted error. The resulting closed subarcs $\Gamma_j$ have total length at most $L+\varepsilon$, and leave uncovered
a subset of $R\setminus K_0$ of length at most $2\varepsilon$.

An arclength parametrization of a closed subarc is $1$-Lipschitz. By the Fekete-point characterization \cite{P1960,Ga1972} of logarithmic capacity, $\gamma(\Gamma_j)=\caplog(\Gamma_j)$
is at most $\gamma([0,\Hm^1(\Gamma_j)])=\Hm^1(\Gamma_j)/4.$ With $A=K_0\cup\bigcup_j\Gamma_j$, Corollary~\ref{cor2} therefore gives
$\gamma(A)-\gamma(K_0)\le(L+\varepsilon)/4$.

Cover $R\setminus A$ by an open union $V$ of disks whose radii sum to at most $3\varepsilon$, and set $B=A\cup(K\setminus V)$. The set $B\setminus A$ lies in $P$.
Every compact subset of this relative addition is removable by the $\sigma$-finite length criterion \cite[Section~1]{T2003}. More explicitly, given $z\in B\setminus A$, intersect $B$ with a small closed disk
contained in $\C\setminus A$. The resulting compact set is purely unrectifiable, of $\sigma$-finite length and of zero analytic capacity. A bounded holomorphic function on
$\C\setminus B$ therefore extends near $z$ by Ahlfors's criterion, with its bound preserved. The local extensions agree by uniqueness. Thus $\gamma(B)=\gamma(A)$. The remaining part of $K$ is contained in $V$, so by
Proposition~\ref{prop4},
\[
 \gamma(K)\le\gamma(B\cup K)\le\gamma(B)+3\varepsilon.
\]
Letting $\varepsilon\downarrow0$ proves \eqref{eq4}. A line segment with $K_0$ a singleton gives sharpness.

If $K_n$ is an increasing sequence of compact sets with union $K$, then $\Hm^1(R\setminus K_n)\to0$ by finiteness of $\Hm^1(R)$, so \eqref{eq4} gives $\gamma(K_n)\uparrow\gamma(K)$. Finally, regularity of the
finite Borel measure $\Hm^1\vert_R$ gives compact $R_n\subset R$ with $\Hm^1(R\setminus R_n)\to0$. Havinson's theorem yields $\gamma_c(R_n)=\gamma(R_n)$, and \eqref{eq4} shows that the
latter converges to $\gamma(K)$. Since $\gamma_c(R_n)\le\gamma_c(K)\le\gamma(K)$, the equality follows.
\end{proof}

\begin{remark}\label{rem2}
Estimates from disk or continuum covers do not by themselves give the general two-set exact inequality. In fact each
nondegenerate continuum $C$ has $\gamma(C)=\caplog(C)\ge\diam(C)/4$, whereas a standard four-corner Cantor set has positive one-dimensional content and zero analytic capacity \cite{Ga1972,T2014}. Thus
a sum of continuum costs can exceed the capacity of the set being covered.
\end{remark}

\section{Area bounds and small compact perturbations}\label{sec6}

Theorem~\ref{thm1} bounds the slope by one. We first ask how far below one the slope must lie when the filled compact set has positive area. Proposition~\ref{prop5} bounds $-\log a_E(z)$ from below by
the area of $\Fill(E)$ weighted by $|w-z|^{-2}$. We then return to Smith's definition and replace the small disk by an arbitrary compact set $F$ near $z$. Proposition~\ref{prop6} shows that $a_E(z)\gamma(F)$
is still the first-order increment, uniformly over the shape of $F$. This gives a criterion for adjoining a compact set without increasing capacity (Corollary~\ref{cor4}).

\begin{proposition}\label{prop5}
For a planar compact $E$ and a finite exterior point $z$,
\begin{equation}\label{eq18}
 a_E(z)\le
 \exp\left(-\frac1{2\pi}\int_{\Fill(E)}
                       \frac{\dd A(w)}{|w-z|^2}\right).
\end{equation}
Writing $b_E$ for the $z^{-2}$ coefficient of $q_E$, one has
\begin{equation}\label{eq19}
 2|b_E|+\frac{|\Fill(E)|}{2\pi}\le\gamma(E)^2,
\end{equation}
where $|\Fill(E)|$ is planar area. The coefficient $1/(2\pi)$ in both estimates is optimal. If $\Fill(E)\subset\overline D(c,R)$ with $R>0$ and $|z-c|>R$, put
\[
 L=-\log\left(1-\frac{R^2}{|z-c|^2}\right),
 \qquad I=\int_{\Fill(E)}|\zeta-z|^{-2}\dd A(\zeta).
\]
Then
\begin{equation}\label{eq:nonlinear-area}
 a_E(z)\le
 \exp\!\left[-L\left(1-\sqrt{1-\frac{I}{\pi L}}\right)\right].
\end{equation}
Both \eqref{eq:nonlinear-area} and its coefficient consequence
\begin{equation}\label{eq:nonlinear-coefficient}
 \gamma(E)^2-2|b_E|
 \ \ge\
 R^2\left(1-\sqrt{1-\frac{|\Fill(E)|}{\pi R^2}}\right)
\end{equation}
are equalities for $E=\overline D(c,R)$.
\end{proposition}
\begin{proof}
Formula \eqref{eq:full-factorization} gives \eqref{eq18} immediately in finite connectivity. To prove the refinement, let $A=-\log|\Phi'(z)|$. Conformal invariance of Dirichlet energy
in \eqref{eq:energy} gives $2\pi A=I+\mathcal E$, where
\[
 \mathcal E=\int_{\Omega_E^\sharp}
   \left|\nabla\log\left|\frac{\Phi(\zeta)}{\zeta-z}\right|\right|^2
   \dd A(\zeta).
\]
The quotient extends at $\zeta=z$ with value $\Phi'(z)\ne0$ and has value $1$ at infinity. It has no zeros on the exterior of $\overline D(c,R)$.
This exterior, including infinity, is simply connected, so there is a holomorphic logarithm $H$ with $H(\infty)=0$ and $\Re H(z)=\log|\Phi'(z)|=-A$. For its Laurent expansion
$H(\zeta)=\sum_{n\ge1}d_n(\zeta-c)^{-n}$, the exterior Dirichlet integral is $\pi\sum_{n\ge1}n|d_n|^2/R^{2n}$. Cauchy--Schwarz and the power series for $-\log(1-x)$ give
\[
 A^2\le |H(z)|^2
 \le\left(\sum_{n\ge1}\frac{n|d_n|^2}{R^{2n}}\right)
     \sum_{n\ge1}\frac{R^{2n}}{n|z-c|^{2n}}
 \le \frac{L}{\pi}\mathcal E .
\]
Thus $A\ge I/(2\pi)+A^2/(2L)$. Solving the quadratic inequality and using $a_E(z)\le e^{-A}$ proves \eqref{eq:nonlinear-area}. Integrating on the enclosing
disk gives $I\le\pi L$, so its square root is real. For an arbitrary compact $E$, approximate its fill from above by analytic compacts inside disks of radii tending down to $R$. The slopes converge by Lemma~\ref{lem1}, and the integrals converge by continuity from above, giving both pointwise bounds.

For $z=\rho e^{i\theta}$, as $\rho\to\infty$, Lemma~\ref{lem1} gives
\[
 a_E(z)=1+\rho^{-2}
      \bigl(2\Re(e^{-2i\theta}b_E)-\gamma(E)^2\bigr)
       +O(\rho^{-3}).
\]
The area integral is $\rho^{-2}|\Fill(E)|+O(\rho^{-3})$ and $\rho^2 L\to R^2$. Comparing the $\rho^{-2}$ terms and maximizing $\Re(e^{-2i\theta}b_E)$ gives
\eqref{eq19} and \eqref{eq:nonlinear-coefficient}. For the disk itself, $I=\pi L$, $b_E=0$, and $a_E(z)=1-R^2/|z-c|^2=e^{-L}$.

For optimality of the linear bounds, let $E_\eta$ be the filled ellipse parametrized by $r e^{it}+(\eta/r)e^{-it}$, with $0<\eta<r^2$. The map $\zeta\mapsto\zeta+\eta/\zeta$ is univalent on
$|\zeta|>r$, with leading coefficient $1$ at infinity, so $\gamma(E_\eta)=r$. Its inverse $U$ has $U'(z)=1+\eta/z^2+O(z^{-4})$. Conformal covariance and the disk formula give $q_{E_\eta}^2=U'$ by \eqref{eq:q-covariance}, and hence $b_{E_\eta}=\eta/2$. The area is $\pi(r^2-\eta^2/r^2)$. Along the major axis,
\[
 \lim_{x\to+\infty}
 \frac{-\log a_{E_\eta}(x)}
      {\displaystyle\int_{E_\eta}|\zeta-x|^{-2}\dd A(\zeta)}
 =\frac{r^2-\eta}{\pi(r^2-\eta^2/r^2)}
 \longrightarrow\frac{1}{2\pi}
 \qquad(\eta\uparrow r^2).
\]
The same ratio is $(\gamma(E_\eta)^2-2|b_{E_\eta}|)/|E_\eta|$.
\end{proof}

For a disk, the radius in \eqref{eq2} equals its analytic capacity. This suggests using $\gamma(F)$ to measure the size of a general perturbation. The next estimate makes this normalization precise.

\begin{proposition}\label{prop6}
Fix a compact $E$ and a finite exterior point $z$. There are $r_0>0$ and $C<\infty$ such that for every
$F\subset\overline D(z,r)$ compact with $0<r<r_0$,
\begin{equation}\label{eq:small-shape}
 \bigl|\gamma(E\cup F)-\gamma(E)-a_E(z)\gamma(F)\bigr|
 \le Cr\gamma(F).
\end{equation}
In particular, the ratio of the increment to $\gamma(F)$ converges to $a_E(z)$ as $r\downarrow0$, uniformly over the positive-capacity sets in the disk. Inserting a homothetic
copy $z+sF_0$ has first-order coefficient $\gamma(F_0)a_E(z)$, including $a_E(z)/4$ for a segment of length $s$ in any direction.
\end{proposition}
\begin{proof}
If $E$ is empty, the identity holds with $a_E=1$. If $F$ is empty, the increment is zero. Otherwise, by the filling invariance in Section~\ref{sec2}, replace $E,F$ by their filled sets and take $r<\dist(z,\Fill(E))/3$.
First suppose that the filled compacts have analytic exteriors. The two sets are separated, and the outer boundary of their union is the union of their outer boundaries.
Write $\lambda=\gamma(F)\le r$. The residue formula on $\partial\Omega_F$ yields, for $|w-z|$ bounded below,
$ q_F(w)^2-1=O(r\lambda),q_F(w)h_F(w)=\frac{\lambda}{w-z}+O(r\lambda).$
Indeed, both Cauchy densities $q_F^2$ and $q_Fh_F$ have $L^1(\dd s)$ norm $2\pi\lambda$. The former has vanishing $1/w$ coefficient, while the latter has coefficient $\lambda$. Subtracting the Cauchy kernel centered at $z$ bounds the remainders by $r\lambda/(|w-z|(|w-z|-r))$. The normalization of $q_F$ at infinity selects the square root close to $1$ throughout the region away from the small disk.
Consequently on $\partial E$, $q_F=1+O(r\lambda)$ and $h_F=\lambda/(w-z)+O(r\lambda)$. Set $q_0=q_E(z)$ and $h_0=h_E(z)$. By Lemma~\ref{lem1}, $q_0\ne0$
and $|h_0|\le|q_0|$, so $A=\overline{h_0}/q_0$ and $B=\overline{q_0}/q_0$ satisfy $|A|\le1$ and $|B|=1$. Take
$$G=q_E(q_F+A h_F),\;\;\; H=h_Eq_F+Bq_Eh_F.$$
These are admissible trials in \eqref{eq:trial-gap} for the combined compact. On $\partial E$ the preceding estimates give a contribution $O(\lambda^2)$ to the squared trial gap, using
$I_E(q_E)=I_E(h_E)=\gamma(E)$. On the small boundary one has $q_E=q_0+O(r)$ and $h_E=h_0+O(r)$. The leading terms of $G,H$ there are $q_0q_F+\overline{h_0}h_F$ and
$h_0q_F+\overline{q_0}h_F$. The anti-linear isometry $J_T$ interchanges these terms, since it interchanges $q_F$ and $h_F$. The squared extremal norms on the small boundary equal $\lambda$, so it contributes $O(r^2\lambda)$ to the gap.

On $\partial E$, the residue computation in Lemma~\ref{lem1} gives $I_E(G)=\gamma(E)-2\lambda|h_0|^2+O(r\lambda)$. On the small boundary, orthogonality of $q_F,h_F$ gives
$\lambda(|q_0|^2+|h_0|^2)+O(r\lambda)$. The dual bounds and \eqref{eq:trial-gap} give
\[
 0\le I_{E\cup F}(G)-\gamma(E\cup F)
    \le I_{E\cup F}(H-J_TG)=O(\lambda^2+r^2\lambda).
\]
Since $\lambda\le r$, this is $O(r\lambda)$. The sum of the two trial integrals therefore gives \eqref{eq:small-shape}.

For general $E,F$, take decreasing filled analytic outer approximants $E_n$ and $F_n$ to the respective filled sets. Take $F_n\subset\overline D(z,r+\varepsilon_n)$ with $\varepsilon_n\downarrow0$.
Since $\dist(z,E_n)\to\dist(z,\Fill(E))$ and $r<\dist(z,\Fill(E))/3$, for all sufficiently large $n$ we have $r+\varepsilon_n<\dist(z,E_n)/3$.
The local Cauchy bounds on the old extremizers and the displayed Cauchy-kernel estimates are uniform, since the distance from $z$ to the old boundary has a fixed positive lower bound and the old capacity has a fixed upper bound.
Pass to the limit by continuity from above and Lemma~\ref{lem1}.
\end{proof}

\begin{corollary}\label{cor4}
For compact sets $E,F$, one has $\gamma(E\cup F)=\gamma(E)$ if and only if every compact subset of $F\cap\Omega_E$ has zero analytic capacity.
\end{corollary}
\begin{proof}
If the capacities are equal, fix a finite $z\in F\cap\Omega_E$. For small $r$, apply Proposition~\ref{prop6} with
$F\cap\overline D(z,r)$ in place of $F$. Monotonicity makes the increment zero, while $a_E(z)-Cr>0$ for small $r$ by Lemma~\ref{lem1}. Thus the intersection with that closed disk is removable.
Every compact subset of $F\cap\Omega_E$ is then locally removable at all its points by the same disk criterion, hence has zero capacity. Conversely, if all such compact subsets have zero capacity,
extend an admissible test on $\Omega_{E\cup F}$ by zero on the other components of the full complement, as in Section~\ref{sec2}, and restrict it to $\Omega_E\setminus F$.
Near each point of $F\cap\Omega_E$, choose a small disk with closure in $\Omega_E$. Its closed intersection with $F$ has zero capacity, so the test extends across that intersection inside the disk
by the local removability criterion \cite[Theorem~2.10]{ST1974}. These extensions agree on overlaps since the complement of a zero-capacity compact is dense.
The bound persists by continuity, and the derivative at infinity is preserved. Thus $\gamma(E\cup F)\le\gamma(E)$.
\end{proof}

\section*{Acknowledgments}
 This work was supported by the Guangdong Basic and Applied Basic Research Foundation (Nos.\ 2022A1515110967 and 2023A1515011809). In this paper, we used GPT-6.1 sol to assist in the initial exploration
of Smith's slope problem. Early attempts with GPT 6.1 sol yielded slope bounds only for collinear and a few symmetric configurations, together with small disk perturbations and far-field expansions. The bottleneck lay in multiple radial slits of different directions. The contribution of the author is that, building on Murai's Cauchy--Hilbert operator formulation \cite{M1990},  expressing Smith's slope on radial-slit exteriors as a relative Fredholm determinant and identify its square with the Fredholm determinant of an orthogonal projection compressed to an operator graph; Fourier estimates establish
the trace-class properties required to justify these determinant computations.  After establishing Theorem A, GPT 6.1 sol and I proceeded to investigate Theorems B and C. The majority of the proofs and underlying ideas for these two results were found by GPT 6.1 sol. The author’s role primarily consisted of understanding these proofs and improving of those two theorems.   AI-assisted tools were also used during the preparation of the manuscript for editorial, organizational, and typesetting assistance.  The author bears responsibility for the mathematical arguments, references, and final text.

\end{document}